\documentclass[12pt]{amsart}

\usepackage[textheight=23.5cm, left=1.2cm, right=1.2cm]{geometry}

\usepackage[english]{babel}
\usepackage[T1]{fontenc}
\usepackage[utf8]{inputenc}

\theoremstyle{plain}
    \newtheorem{theorem}{Theorem}[section]
    \newtheorem{lemma}[theorem]{Lemma}

\theoremstyle{definition}

\usepackage[
    draft = false,
    unicode = true,
    colorlinks = true,
    allcolors = blue,
    hyperfootnotes = true,
    citecolor = red
]{hyperref}
\usepackage{amsfonts,amsmath,amssymb,amscd,amsthm,url}
\usepackage[inline]{enumitem}
\usepackage{graphicx,epsf,afterpage}
\usepackage{soul}
\usepackage{multicol}
\usepackage{array}
\usepackage{braket}
\usepackage{epigraph}
\usepackage{rotating}

\usepackage[
backend=biber,
style=numeric-comp,
sorting=none,
giveninits=true
]{biblatex}
\usepackage{csquotes}

\usepackage{xcolor}

\usepackage{ulem}

\usepackage{tikz}
\usetikzlibrary{positioning, calc, intersections, through, arrows, matrix, chains, math}
\usetikzlibrary{shadows}
\usetikzlibrary{graphs}
\usetikzlibrary{graphs.standard}
\usetikzlibrary{patterns}
\usetikzlibrary{decorations.pathreplacing,calligraphy,backgrounds}

\newcommand\norm[1]{\ensuremath{\left\lVert#1\right\rVert}}
\newcommand\abs[1]{\ensuremath{\left\lvert#1\right\rvert}}

\renewcommand{\Pr}{\mathrm{P}}

\DeclareMathOperator{\Expect}{\mathbb{E}}

\newcommand{\Lcal}{\mathcal{L}}
\newcommand{\Hcal}{\mathcal{H}}

\newcommand{\Fcal}{\mathcal{F}}

\newcommand{\R}{\ensuremath{\mathbb{R}}}

\newcommand{\N}{\ensuremath{\mathbb{N}}}

\renewcommand{\leq}{\leqslant}

\newcounter{mcnt}

\newcounter{wordcnt}

\begin{document}

\title[The CLT for random exponents on a Hilbert space in WOT]{The Central Limit Theorem \\ for random exponents on a Hilbert space \\ in the Weak Operator Topology}

\author{S.\,V. Dzhenzher}

\begin{abstract}
     We consider random linear continuous operators $\Omega \to \mathcal{L}(\mathcal{H}, \mathcal{H})$ on a Hilbert space $\mathcal{H}$.
     For example, such random operators may be random quantum channels.
     The Central Limit Theorem is known for the sums of i.i.d. random operators.
     Instead of the sum, there may be considered the composition of random exponents $e^{A_i/n}$.
     We obtain the Central Limit Theorem in the Weak Operator Topology for centralized and normalized random exponents of i.i.d. linear continuous operators on a Hilbert space.
\end{abstract}

\thanks{
Moscow Institute of Physics and Technology 141701, Dolgoprudny, Russia. sdjenjer@yandex.ru \\
I am grateful to V.\,Zh.\,Sakbaev for many helpful discussions and suggestions.
}

\maketitle

\thispagestyle{empty}


\section{Introduction}

The theory of compositions of random linear operators is constructed and developed in a lot of publications \cite{Meta, Furstenberg, Tutubalin}. 
A random quantum dynamics can be considered either as a composition of random linear operators in a Hilbert space \cite{Kempe, OSS19}, or as a composition of random quantum channels \cite{Aaron, Hol, GOSS-2022, New-23, N-23}.
See in \cite{Berezin, OSS-2016, Orlov2, Aaron, Hol, Pechen-Volovich, GOSS21, GOSS-2022, Joye, Pechen} how random operator-valued processes arise.

The Law of Large Numbers (LLN) and the Central Limit Theorem (CLT) for random variables acting in a Banach space is well studied \cite{Ledoux-Talagrand-1991-probability, Hoffman-Jorgensen-Pisier}.
Distributions of products of random linear operators have been studied in \cite{Oseledets, Tutubalin, Skorokhod, New-23}.
Different versions of limit theorems of LLN and CLT type for compositions of i.i.d. random operators were obtained in \cite{Berger, Watkins-1984, Watkins-1985, GOSS-2022, DzhenzherSakbaev24, DzhenzherSakbaev25, DzhenzherSakbaev25-SLLN, Orlov-Sakbaev-Shmidt-2021, Sakbaev-Shmidt-Shmidt-2022, SSSh23}.

Here are some versions and approaches for the CLT-type problems.
\begin{itemize}[nosep]
    \item 
    In \cite{Beigi2025} the convolutions of quantum states are considered. It is shown that the speed of convergence of these convolutions to a Gaussian state by the Schatten $1$-norm is \(O(1/\sqrt{n})\).
    
    \item In \cite{Orlov-Sakbaev-Shmidt-2021} it is shown that the weak convergence of measures is equivalent to the pointwise convergence of convolution operators.
    Different topologies of test functions instead of continuous bounded functions are considered.
    In those topologies the CLT-type results are obtained for the random operator-valued processes.
    
    \item In \cite{Berger} there are considered the random matrices of the kind \(I + \frac{A_i}{\sqrt{n}} + \frac{B_i}{n} + O(n^{-3/2})\).
    It is shown there that their products converge uniformly (in distribution) to the log-normal random distribution given by some stochastic equation.
    
    \item In \cite{Watkins-1984, Watkins-1985} the composition of the kind \(e^{A(\xi_1)/n}\ldots e^{A(\xi_{\lfloor n^2t\rfloor})/n}\) is considered.
    The weak convergence of these processes is proved.
\end{itemize}

In \cite{DzhenzherSakbaev25-SLLN} were studied the compositions of random semigroups of the form \(e^{A_it/n}\) for operators on Banach spaces.
The Strong Law of Large Numbers was obtained for them in the spirit of Chernoff's Theorem \cite{Chernoff-1968}; it was obtained in the Weak Operator Topology for all Banach spaces, and in the Strong Operator Topology for uniformly smooth Banach spaces.
Following \cite{Orlov-Sakbaev-Shmidt-2021} it is only natural to centralize the random compositions in \cite{DzhenzherSakbaev25-SLLN}, multiply them by \(\sqrt{n}\), and obtain the CLT for them as is done in classical Probability.
In this paper, we obtain such CLT for operators on a Hilbert space.
Precisely, we obtain that the sequence
\[
    \sqrt{n}\left(e^{A_1/n} \ldots e^{A_n/n} - e^{\Expect A}\right)
\]
of operators converges in distribution in the Weak Operator Topology to a random Gaussian operator.

The structure of the article is the following.
In~\S\ref{s:main}, the main result (Theorem~\ref{t:clt}) is given.
In~\S\ref{s:proof}, the proof of Theorem~\ref{t:clt} is given.

\section{Background and statements}\label{s:main}

Let \((\Omega, \Fcal, \Pr)\) be a probability space with a complete $\sigma$-additive measure $\Pr$.

For a Hilbert space \(\Hcal\) denote by \(\Lcal(\Hcal)\) the space of bounded linear operators on \(\Hcal\).
The map \(A:\Omega\to\Lcal(\Hcal)\) is called a \textbf{\emph{random operator}} if it is \emph{Bochner measurable}, i.e. if there is a sequence of finitely-valued maps \(\Omega\to\Lcal(\Hcal)\) converging to \(A\) uniformly almost surely.
The \textbf{\emph{expected value}} \(\Expect A\) of a random operator is its \emph{Bochner integral}.
It is well-known \cite{DiestelUhl, Yosida1995} that \(\Expect A\) exists if and only if the random variable \(\norm{A}\) is Lebesgue integrable.

The independence of random operators is defined analogously to the independence of random variables.
See \cite{DzhenzherSakbaev25-SLLN} for the proof that the composition of random operators is a random operator, and that the expected value of the composition of independent random operators is the composition of their expected values.

For a linear bounded positive symmetric operator \(\Sigma\) on \(\Hcal^{\otimes2}\) we denote by \(N(0,\Sigma)\) the distribution of a zero-mean random Gaussian operator $G$ with the covariance \(\Sigma\); that is, for any \(x,y\in\Hcal\) the random variable \(\braket{y, Gx}\) is a zero-mean Gaussian variable with the variance equal to \(\Braket{y^{\otimes2},\Sigma x^{\otimes2}}\).

\begin{theorem}\label{t:clt}
    Let \(\Hcal\) be a Hilbert space over \(\R\).
    Let \(A\colon\Omega\to\Lcal(\Hcal)\) be such that $\norm{A}$ is bounded.
    Let $A,A_1, A_2, \ldots$ be a sequence of i.i.d. random operators.
    Then the sequence
    \[
        \sqrt{n}\left(e^{A_1/n} \ldots e^{A_n/n} - e^{\Expect A}\right)
    \]
    of random operators converges in distribution in WOT to \(N(0,\Sigma)\) as \(n\to\infty\), where
    \[
        \Sigma := \int_{s=0}^1 \left(e^{\Expect As}\right)^{\otimes2} \Expect (A-\Expect A)^{\otimes2} \left(e^{\Expect A(1-s)}\right)^{\otimes2}\,ds
    \]
    is the Bochner integral.
    In other words, for any \(x,y\in \Hcal\) we have
    \[
        \Braket{y,\sqrt{n}\left(e^{A_1/n} \ldots e^{A_n/n} - e^{\Expect A}\right)x} \xrightarrow[n\to\infty]{d} N\left(0,\Braket{y^{\otimes2},\Sigma x^{\otimes2}}\right).
    \]
\end{theorem}

Note that in \cite[Theorem~2.5]{DzhenzherSakbaev25-SLLN} it is shown that \(\sup\limits_{t\in[\,0,T\,]}\abs{\Braket{y,\left(e^{A_1t/n} \ldots e^{A_nt/n} - e^{\Expect At}\right)x}}\) converges almost surely to \(0\).
Therefore, it would be interesting to obtain the convergence in distribution as in Theorem~\ref{t:clt}, but for the \(\sup\) for semigroups instead of just exponents, or as in \cite{Berger}.

\section{Proof of the main result}\label{s:proof}

\begin{lemma}[{cf.~\cite[Theorem~3]{Galkin-Remizov-2025}}]\label{l:speed}
    Let \(\Hcal\) be a Hilbert space.
    Let \(A\colon\Omega\to\Lcal(\Hcal)\) such that $\norm{A}$ is bounded.
    Then for any \(n \in \N\) and \(1 \leq k \leq n\)
    \[
        \norm{\left(\Expect e^{A/n}\right)^k - e^{\Expect Ak/n}} = O(1/n)
        \quad\text{as \(n\to\infty\)},
    \]
    with the constant in $O(1/n)$ independent of $k = k(n)$.
\end{lemma}

\begin{proof}
    By the Taylor's formula for Banach spaces, both \(\Expect e^{A/n}\) and \(e^{\Expect A/n}\) could be expressed as
    \(
        I + \frac{\Expect A}{n} + O\left(1/n^2\right).
    \)
    Then
    \[
        \norm{\Expect e^{A/n} - e^{\Expect A/n}} = O\left(1/n^2\right).
    \]
    Now the result follows from the formula
    \[
        X^k - Y^k = \sum_{j=0}^{k-1} X^{k-j-1}(X-Y)Y^j
    \]
    for \(X := \Expect e^{A/n}\) and \(Y := e^{\Expect A/n}\).
\end{proof}

\begin{proof}[Proof of Theorem~\ref{t:clt}]
    Fix arbitrary \(y,x\in\Hcal\).
    Denote
    \[
        \xi_n := \sqrt{n}(e^{A_1/n} \ldots e^{A_n/n} - e^{\Expect A}) \quad
        \text{and} \quad
        S_n := \frac{1}{\sqrt{n}}\sum_{k=1}^n e^{\Expect A\frac{k-1}{n}} (A_k - \Expect A) e^{\Expect A\frac{n-k}{n}}.
    \]
    We want to prove that \(\norm{\xi_nx - S_nx}\xrightarrow[n\to\infty]{\Pr}0\), and that $S_n$ is a triangular array of martingales w.r.t \(\Fcal_{n,k} = \Fcal_k = \sigma(A_1,\ldots, A_k)\) for which the following conditions hold:
    \begin{gather}
        \label{eq:cov-conv}
        \tag{Cov}
        \sum_{k=1}^n \Expect^{\Fcal_{k-1}} d_{n,k}\otimes d_{n,k} \xrightarrow[n\to\infty]{\Pr} \Sigma, \\
        \label{eq:lind}
        \tag{Lind}
        \forall\varepsilon>0\quad\sum_{k=1}^n \Expect^{\Fcal_{k-1}} \norm{d_{n,k}}^2 \mathcal{I}\{\norm{d_{n,k}} > \varepsilon\} \xrightarrow[n\to\infty]{\Pr} 0,
    \end{gather}
    where \(d_{n,k}\) are the summands in $S_n$.
    Then the result will follow by the CLT \cite[Corollary~3.1]{Hall-Heyde} for real-valued martingales.

    Let us show that $S_n$ is a triangular array of martingale differences.
    Indeed, \(\Expect^{\Fcal_{k-1}} d_{n,k} = 0\) since $A_k$ is independent of \(\Fcal_{k-1}\).

    Let us show that condition~\eqref{eq:cov-conv} holds for $S_n$.
    Since \(A_k\) is independent of \(\Fcal_{k-1}\), we have
    \[
        \Expect^{\Fcal_{k-1}} d_{n,k}^{\otimes2} = \Expect d_{n,k}^{\otimes2} = \frac{1}{n}
        \left(e^{\Expect A\frac{k-1}{n}}\right)^{\otimes2} \Expect (A-\Expect A)^{\otimes2} \left(e^{\Expect A\frac{n-k}{n}}\right)^{\otimes2}.
    \]
    As $n \to \infty$, $k/n \to s$, $\frac{1}{n} \to ds$, we have that the sums of these expectations converge to \(\Sigma\) and~\eqref{eq:cov-conv} holds.
    
    Let us show that conditional Lindeberg condition~\eqref{eq:lind} holds for $S_n$.
    Let \(\rho>0\) be the constant such that \(\norm{A}\leqslant \rho\).
    Hence \(\norm{e^{A}} \leq e^{\rho}\), and
    \[
        \norm{d_{n,k}} \leq \frac{2\rho}{\sqrt{n}} e^{\rho \frac{n-1}{n}}.
    \]
    Thus for any \(\varepsilon > 0\) for large enough $n$ we have \(\{\norm{d_{n,k}} > \varepsilon\} = \varnothing\), and~\eqref{eq:lind} holds.

    It remains to show that \(\norm{\xi_nx - S_nx} \xrightarrow[n\to\infty]{\Pr} 0\).
    Since by Lemma~\ref{l:speed} for \(k=n\) we have that
    \[
        \sqrt{n}\norm{(\Expect e^{A_1/n} \ldots e^{A_n/n} - e^{\Expect A})x} \xrightarrow[n\to\infty]{} 0,
    \]
    it is sufficient to show that
    \[
        \norm{S_nx -\xi_n'x} \xrightarrow[n\to\infty]{\Pr} 0,
    \]
    where
    \[
        \xi_n' := \sqrt{n}\left(e^{A_1/n} \ldots e^{A_n/n} - \left(\Expect e^{A/n}\right)^n\right) =
        \sqrt{n}\sum_{k=1}^ne^{A_1/n} \ldots e^{A_{k-1}/n}\left(e^{A_k/n} - \Expect e^{A/n}\right) \left(\Expect e^{A/n}\right)^{n-k}.
    \]
    Applying the Taylor's formula for \(e^{A_k/n} - \Expect e^{A/n}\) we have
    \[
        \xi_n' = S_n' + R_n,
    \]
    where
    \[
        S_n' := \frac{1}{\sqrt{n}} \sum_{k=1}^ne^{A_1/n} \ldots e^{A_{k-1}/n}\left(A_k - \Expect A\right) \left(\Expect e^{A/n}\right)^{n-k}
    \]
    and \(\norm{R_n} = \sqrt{n} \cdot n \cdot e^{\rho} \cdot O(1/n^2) = O(1/\sqrt{n})\).
    So it remains to show that \(\norm{S_nx - S_n'x} \xrightarrow[n\to\infty]{\Pr} 0\).
    By the Markov's inequality it is sufficient to show that
    \[
        \Expect \norm{S_nx - S_n'x}^2 \xrightarrow[n\to\infty]{} 0.
    \]
    Denote the summands in $S_n'$ by \(d_{n,k}'\).
    Now we have that
    \[
        \Expect \norm{S_nx - S_n'x}^2 = \sum_{k,\ell} \Expect \braket{d_{n,k}x - d_{n,k}'x, d_{n,\ell}x - d_{n,\ell}'x} = 
        \sum_{k=1}^n \Expect \norm{d_{n,k}x - d_{n,k}'x}^2,
    \]
    where in the last equality we used that for \(k > \ell\) by the law of total probability
    \[
        \Expect \braket{d_{n,k}x - d_{n,k}'x, d_{n,\ell}x - d_{n,\ell}'x} =
        \Expect \Braket{\Expect^{\Fcal_{k-1}}(d_{n,k}x - d_{n,k}'x), d_{n,\ell}x - d_{n,\ell}'x} =0,
    \]
    since \(\Expect^{\Fcal_{k-1}} d_{n,k} = \Expect^{\Fcal_{k-1}} d_{n,k}' = 0\).
    Finally, in order to complete the proof it is sufficient to show that
    \[
        \Expect \norm{d_{n,k}x - d_{n,k}'x}^2 = O\left(1/n^2\right),
    \]
    where the constant is independent of $k=k(n)$.
    
    Denote
    \[
        M_k := e^{A_1/n} \ldots e^{A_k/n}x - \left(\Expect e^{A/n}\right)^kx.
    \]
    Notice that up to unimportant constant depending on $\rho$ we have
    \[
        \sqrt{n}\norm{d_{n,k}x - d_{n,k}'x} \leq \norm{M_{k-1}} +
        \norm{\left(\Expect e^{A/n}\right)^{k-1}x - e^{\Expect A\frac{k-1}{n}}x} +
        \norm{\left(\Expect e^{A/n}\right)^{n-k}x - e^{\Expect A\frac{n-k}{n}}x}.
    \]
    Here the last two summands are \(O(1/n)\) by Lemma~\ref{l:speed}.
    Then
    \[
        n\norm{d_{n,k}x - d_{n,k}'x}^2 \leq \norm{M_{k-1}}^2 + 
        \norm{M_{k-1}} \cdot O(1/n) + O(1/n^2).
    \]
    Analogously to \cite{DzhenzherSakbaev25-SLLN}, for integer $0 \leq m \leq k$ denote \([k]:= \{1,\ldots,k\}\), and denote by \(\binom{[k]}{m}\) the family of $m$-element subsets of \([k]\).
    For \(\{i_1 < \ldots < i_m\} \in \binom{[k]}{m}\) denote
    \[
        F_{k,\{i_1, \ldots, i_m\}} := (\Expect e^{A/n})^{i_1-1} \left(e^{A_{i_1}/n} - \Expect e^{A/n}\right) (\Expect e^{A/n})^{i_2-i_1-1} \ldots \left(e^{A_{i_m}/n} - \Expect e^{A/n}\right) (\Expect e^{A/n})^{k-i_m}.
    \]
    Denoting
    \[
        D_{k,m} := \sum_{\substack{P \subset [k] \\ \max P = m}} F_{k,P}x,
    \]
    analogously to \cite{DzhenzherSakbaev25-SLLN} we obtain that
    \[
        M_k = \sum_{m=1}^k D_{k,m}
    \]
    is the martingale w.r.t. to \(\Fcal_k = \sigma(A_1,\ldots,A_k)\).
    By \cite[Lemma~3.3]{DzhenzherSakbaev25-SLLN} we have
    \[
        \norm{F_{k,P}} \leq \left(\frac{2\rho}{n}\right)^{\abs{P}} e^{k\rho/n}.
    \]
    Hence up to $\norm{x}$
    \[
        \norm{D_{k,m}} \leq
        \sum_{j=1}^m \binom{m-1}{j-1} \left(\frac{2\rho}{n}\right)^j e^{\rho} =
        \frac{2\rho}{n} e^{\rho} \left(1+\frac{2\rho}{n}\right)^{m-1} \leq \frac{2\rho}{n}e^{3\rho}.
    \]
    Since $\Hcal$ is $2$-smooth, we may apply \cite[Theorem~4.52]{Pisier-2016-martingales}; see \cite[Burkholder-type Theorem~4.1]{DzhenzherSakbaev25-SLLN} for the exposition adapted to this paper.
    In any case, up to some constant independent of $k$ we obtain
    \[
        \Expect\norm{M_k} \leq \sqrt{n \cdot \left(\frac{2\rho}{n}e^{3\rho}\right)^2} = O(1/\sqrt{n})
    \quad\text{and}\quad
        \Expect\norm{M_k}^2 \leq n \cdot \left(\frac{2\rho}{n}e^{3\rho}\right)^2 = O(1/n).
    \]
    Thus
    \[
        \norm{d_{n,k}x - d_{n,k}'x}^2 = O\left(1/n^2\right),
    \]
    which completes the proof.
\end{proof}

\printbibliography

@article{Berezin,
  author    = {F.~A. Berezin},
  title     = {Non-Wiener functional integrals},
  journal   = {Theoret. and Math. Phys.},
  volume    = {6},
  number    = {2},
  pages     = {141--155},
  year      = {1971}
}

@article{Berger,   
author       ={M.~A. Berger },   
title        ={Central limit theorem for products of random matrices }, 
journal      ={Trans.AMS.},   
volume       ={285},
series       ={2},
pages        ={777--803 },   
year         ={1984 }   
}

@Article{Beigi2025,
    author={Beigi, Salman
    and Mehrabi, Hami},
    title={Toward Optimal Convergence Rates for the Quantum Central Limit Theorem},
    journal={Annales Henri Poincar{\'e}},
    year={2025},
    month={Jul},
    day={31},
    issn={1424-0661},
    doi={https://doi.org/10.1007/s00023-025-01609-4}
}

@article{Chernoff-1968,
title = {Note on product formulas for operator semigroups},
journal = {Journal of Functional Analysis},
volume = {2},
number = {2},
pages = {238-242},
year = {1968},
issn = {0022-1236},
doi = {https://doi.org/10.1016/0022-1236(68)90020-7},
author = {Paul R Chernoff}
}

@article{DzhenzherSakbaev24,
    author = {S.~V. Dzhenzher and V.~Zh. Sakbaev},
    title = {Quantum law of large numbers for {B}anach spaces},
    journal = {Lobachevskii Journal of Mathematics},
    volume={45},
    pages={2485--2494},
    doi = {https://doi.org/10.1134/S1995080224603114},
    year = {2024}
}

@article{DzhenzherSakbaev25,
    author = {S.~V. Dzhenzher and V.~Zh. Sakbaev},
    title = {The law of large numbers for discrete generalized quantum channels},
    journal = {Lobachevskii Journal of Mathematics (accepted)},
    url = {arXiv:2504.10033},
    year = {2025}
}

@article{DzhenzherSakbaev25-SLLN,
    author = {S.~V. Dzhenzher and V.~Zh. Sakbaev},
    title = {The Strong Law of Large Numbers for random semigroups on uniformly smooth Banach spaces},
    url = {arXiv:2507.07658}
}

@article{Aaron,
  author = {Aharonov, Y. and Davidovich, L. and Zagury, N. },
  title = {Quantum random walks},
  journal = {Physical Review A},
  volume = {48},
  number = {2},
  pages = {1687--1690},
  year = {1993}
}

@incollection{GOSS21,
  author = {Gough, J. and Orlov, Yu.~N. and Sakbaev, V.~Zh. and Smolyanov, O.~G.},
  title = {Random quantization of Hamiltonian systems},
  booktitle = {Doklady Mathematics},
  volume = {103},
  pages = {122--126},
  year = {2021},
  publisher = {Springer}
}

@article{Hoffman-Jorgensen-Pisier,
author = {J. Hoffmann-Jorgensen and G. Pisier},
title = {{The Law of Large Numbers and the Central Limit Theorem in Banach Spaces}},
volume = {4},
journal = {The Annals of Probability},
number = {4},
publisher = {Institute of Mathematical Statistics},
pages = {587--599},
year = {1976},
doi = {https://doi.org/10.1214/aop/1176996029}
}

@book{Meta,
  author = {M.~L. Mehta},
  title = {Random matrices},
  publisher = {Elsevier},
  address = {Amsterdam},
  year = {2004}
}

@article{Pechen,
  author = {A.~N. Pechen},
  title = {Quantum stochastic equation for a test particle interacting with a dilute Bose gas},
  journal = {Journal of Mathematical Physics},
  volume = {45},
  number = {1},
  pages = {400--417},
  year = {2004}
}

@misc{N-23,
  author = {S. Dhamapurkar and O. Dahlsten},
  title = {Quantum walks as thermalizations, with application to fullerene graphs},
  howpublished = {arXiv:2304.01572},
  year = {2023}
}

@article{Furstenberg,
  author = {H. Furstenberg},
  title = {Non-commuting random products},
  journal = {Transactions of the American Mathematical Society},
  volume = {108},
  number = {3},
  pages = {377--428},
  year = {1963}
}

@article{Joye,
  author = {A. Joye},
  title = {Random time-dependent quantum walks},
  journal = {Communications in Mathematical Physics},
  volume = {307},
  number = {1},
  pages = {65--100},
  year = {2011}
}

@article{Pechen-Volovich,
  author = {A.~N. Pechen and I.~V. Volovich},
  title = {Quantum multipole noise and generalized quantum stochastic equations},
  journal = {Infinite Dimensional Analysis, Quantum Probability and Related Topics},
  volume = {5},
  number = {4},
  pages = {441--464},
  year = {2002}
}

@article{Galkin-Remizov-2025,
author = {Oleg E. Galkin and Ivan D. Remizov},
year = {2025},
month = {02},
pages = {929--943},
title = {Upper and lower estimates for rate of convergence in the Chernoff product formula for semigroups of operators},
volume = {265},
issue={2},
journal = {Israel Journal of Mathematics},
doi = {https://doi.org/10.1007/s11856-024-2678-x}
}

@article{OSS-2016,
  author    = {Yu.~N. Orlov and V.~Zh. Sakbaev and O.~G. Smolyanov},
  title     = {Unbounded random operators and Feynman formulae},
  journal   = {Izv. Math.},
  volume    = {80},
  number    = {6},
  pages     = {1131--1158},
  year      = {2016}
}

@article{Orlov2,
  author    = {Yu.~N. Orlov},
  title     = {Evolution equation for Wigner function for linear quantization},
  journal   = {Keldysh Institute preprints},
  volume    = {040},
  pages     = {22},
  year      = {2020}
}

@article{GOSS-2022,
  author    = {J.~Gough and Yu.~N. Orlov and V.~Zh. Sakbaev and O.~G. Smolyanov},
  title     = {Markov Approximations of the Evolution of Quantum Systems},
  journal   = {Dokl. Math.},
  volume    = {105},
  number    = {2},
  pages     = {92--96},
  year      = {2022}
}

@misc{New-23,
  author       = {L. Pathirana and J. Schenker},
  title        = {Law of large numbers and central limit theorem for ergodic quantum processes},
  howpublished = {arXiv:2303.08992},
  year         = {2023}
}

@article{Skorokhod,
  author    = {A.~V. Skorokhod},
  title     = {Products of independent random operators},
  journal   = {Russian Math. Surveys,},
  volume    = {38},
  number    = {(4)},
  pages     = {291--318},
  year      = {1983}
}

@article{Tutubalin,
   author={V.~N. Tutubalin}, 
   title={A local limit theorem for products of random matrices}, 
   journal={Probab. Theory and Appl.,}, 
   volume={22}, 
   number={2}, 
   pages={203--214}, 
   year={1978} 
}

@article{Oseledets,
   author={V.~I. Oseledets}, 
   title={Markov chains, skew products and ergodic theorems for general dynamic systems}, 
   journal={Probab. Theory and Appl.,}, 
   volume={10}, 
   number={3}, 
   pages={551--557}, 
   year={1965} 
}

@article{OSS19,
   author={Yu.~N. Orlov and V.~Zh. Sakbaev and O.~G. Smolyanov}, 
   title={Feynman Formulas and the Law of Large Numbers for Random One-Parameter Semigroups}, 
   journal={Proc.Steklov Inst.Math.}, 
   volume={306}, 
   pages={196--211}, 
   year={2019} 
}

@article{SSSh23,   
author       ={Yu.~N. Orlov and V.~Zh. Sakbaev and E.~V. Shmidt },   
title        ={Compositions of Random Processes in a Hilbert Space and Its Limit Distribution },   
journal      ={Lobachevskii J.Math.},
volume       ={44},
series       = {4},
pages        ={1432--1447 },   
year         ={2023 }   
}

@article{Sakbaev-Shmidt-Shmidt-2022,
author={Sakbaev, V. Zh.
and Shmidt, E. V.
and Shmidt, V.},
title={Limit Distribution for Compositions of Random Operators},
journal={Lobachevskii Journal of Mathematics},
year={2022},
month={Jul},
day={01},
volume={43},
number={7},
pages={1740-1754},
issn={1818-9962},
doi={https://doi.org/10.1134/S199508022210033X}
}

@article{Orlov-Sakbaev-Shmidt-2021,
author={Orlov, Yu. N.
and Sakbaev, V. Zh.
and Shmidt, E. V.},
title={Operator Approach to Weak Convergence of Measures and Limit Theorems for Random Operators},
journal={Lobachevskii Journal of Mathematics},
year={2021},
month={Oct},
day={01},
volume={42},
number={10},
pages={2413-2426},
issn={1818-9962},
doi={https://doi.org/10.1134/S1995080221100188}
}

@book{Hol,
author       ={A.~S. Holevo },  
title        ={Quantum probability and quantum statistics },  
publisher    ={Probability theory -8, Itogi Nauki i Tekhniki.Series Sovrem.Probl.Mat.Fund.Napr.},    
volume       ={83} ,    
address      ={VINITI,Moscow },    
pages        ={5--132 },    
year         ={1991 }    
}

@article{Kempe,
author       ={J. Kempe },  
title        ={Quantum random walks: an introductory overview },  
journal      ={Contemp.Phys.},    
volume       ={44(4)},    
pages        ={307--327},    
year         ={2003 }    
}

@article{Watkins-1984,
 author = {Watkins, Joseph C.},
 title = {A central limit problem in random evolutions},
 fjournal = {The Annals of Probability},
 journal = {Ann. Probab.},
 issn = {0091-1798},
 volume = {12},
 pages = {480--513},
 year = {1984},
 language = {English},
 doi = {https://doi.org/10.1214/aop/1176993302},
}

@article{Watkins-1985,
title = {Limit theorems for stationary random evolutions},
journal = {Stochastic Processes and their Applications},
volume = {19},
number = {2},
pages = {189-224},
year = {1985},
issn = {0304-4149},
doi = {https://doi.org/10.1016/0304-4149(85)90025-0},
author = {Joseph C. Watkins},
}

@Inbook{Yosida1995,
author="Yosida, K{\^o}saku",
title="Strong Convergence and Weak Convergence",
bookTitle="Functional Analysis",
year="1995",
publisher="Springer Berlin Heidelberg",
address="Berlin, Heidelberg",
pages="119--145",
isbn="978-3-642-61859-8",
doi="https://doi.org/10.1007/978-3-642-61859-8_6"
}

@incollection{Hall-Heyde,
title = {3 - The Central Limit Theorem},
editor = {P. Hall and C.C. Heyde},
booktitle = {Martingale Limit Theory and its Application},
publisher = {Academic Press},
pages = {51-96},
year = {1980},
series = {Probability and Mathematical Statistics: A Series of Monographs and Textbooks},
issn = {00795607},
doi = {https://doi.org/10.1016/B978-0-12-319350-6.50009-8},
url = {https://www.sciencedirect.com/science/article/pii/B9780123193506500098},
author = {P. Hall and C.C. Heyde}
}

@book{DiestelUhl,
  title={Vector Measures},
  author={Diestel, J. and Uhl, J.~J.},
  isbn={9780821815151},
  lccn={lc77009625},
  series={Mathematical surveys and monographs},
  url={https://books.google.ru/books?id=fEKZAwAAQBAJ},
  year={1977},
  publisher={American Mathematical Society}
}

@book{Ledoux-Talagrand-1991-probability,
  title={Probability in Banach Spaces: Isoperimetry and Processes},
  author={Ledoux, M. and Talagrand, M.},
  isbn={9783540520139},
  lccn={lc91010765},
  series={A Series of Modern Surveys in Mathematics Series},
  url={https://books.google.ru/books?id=cyKYDfvxRjsC},
  year={1991},
  publisher={Springer}
}

@book{Pisier-2016-martingales,
  title={Martingales in Banach Spaces},
  author={Pisier, G.},
  isbn={9781107137240},
  lccn={2016427432},
  series={Cambridge Studies in Advanced Mathematics},
  url={https://books.google.ru/books?id=vJErDAAAQBAJ},
  year={2016},
  publisher={Cambridge University Press}
}

\end{document}